\documentclass[12pt]{article}
\usepackage{a4wide}
\usepackage{latexsym, amsmath, amsfonts, amsthm, amssymb}

\theoremstyle{plain}

\newtheorem*{theo*}{Theorem}

\newcommand{\R}{\ensuremath{\mathbb{R}}} 

\newcommand{\Var}{{\rm Var}} 

\def\id{\textrm{Id}} 
\def\Hess{\mathop{\rm Hess}\nolimits} 
\def\tr{\mathop{\rm tr}\nolimits} 

\def\ric{\mathop{\rm Ric}\nolimits}
\def\vol{\mathop{\rm vol}\nolimits}

\def\eps{\varepsilon}
\def\te{\theta}
\newcommand{\upchi}{\raise1pt\hbox{$\chi$}}

\newcommand{\be}{\begin{equation}}
\newcommand{\ee}{\end{equation}}
\def\benu{\begin{enumerate}}
\def\eenu{\end{enumerate}}
\newcommand{\AND}{\qquad {\rm and} \qquad}

\title{A transport inequality on the sphere \\ obtained by mass transport}
\author{Dario Cordero\,-Erausquin 
}
\date{}      

\begin{document}

\maketitle
\begin{abstract}
Using McCann's transportation map, we establish a transport inequality on  compact manifolds with positive Ricci curvature. This inequality contains the sharp  spectral comparison estimates. 
\end{abstract}

\section{Introduction}

Extending the mass transportation approach to sharp Sobolev type inequalities from Euclidean space to curved geometries remains a challenging problem. In the present note, we propose a new twist in the classical transportation technique that allows for a transport inequality which contains sharp Poincar\'e inequalities.

The method applies to a (compact) Riemannian manifold  of dimension $n\ge 2$ having a lower bound on the Ricci curvature of the form $\ric \ge (n-1)k^2 {\rm g}$ with $k>0$ and ${\rm g}$ the Riemannian metric. By scaling the distances, we can always assume that $k=1$. 

So, in the rest of the paper $M=(M, {\rm g})$ will stand for an $n$-dimensional  Riemannian manifold satisfying 
\be\label{riccibound}
\ric \ge (n-1)\, {\rm g}.
\ee
The  main example is the usual sphere $S^n\subset \R^{n+1}$. The interest, perhaps,   in stating a result under the condition~\eqref{riccibound}, even if one aims at the sphere only, is that it makes it clear that we will not use any of the algebraic properties of the sphere. Our computations are modeled on the sphere case;  the extension to the   situation given by~\eqref{riccibound} relies on Bishop comparison's estimates only.
We will denote by  $d\sigma = d\vol / \vol(M) $ the Riemannian volume measure normalized to be a probability measure. The distance will be denoted $d$; recall that $M$ has diameter smaller than $\pi$.

A simple but important result is that, on such manifold $M$, the  spectral gap for the Laplacian satisfies $\lambda_1 \ge n$. Equivalently, one has the following  Wirtinger-Poincar\'e inequality: for every Lipschitz function $g$ on $M$,
\be\label{poinc}
\Var_{\sigma} (g) := \int \Big(g -\int g\, d\sigma\Big)^2 \, d\sigma \le \frac{1}{n} \int |\nabla g|^2 \, d\sigma.
\ee

The $L^2$ proof of this inequality as done by Lichnerovicz using Bochner's formula  is rather short and elementary. In the particular case of the sphere, one can also use the expansion of $g$ in the spherical harmonics basis; moreover, in this case, equality holds for linear functions, which are eigenfunctions for the spherical Laplacian. 

It is well known that Poincar\'e inequalities are not well suited to mass transport techniques. However, in the Euclidean case and under appropriate curvature assumptions, one can prove very easily using mass transport (Brenier map) \emph{stronger} inequalities such as \emph{transport inequalities} or \emph{logarithmic Sobolev inequalities} (see~\cite{C02}).
So it is quite annoying that no mass transport proof of the sharp log-Sobolev inequality (see~\cite{Le}), say, is available on $M$. Indeed, the straightforward adaptation of the techniques from Euclidean space leads to a log-Sobolev inequality with a constant $(n-1)$ in place of the expected constant $n$. Similarly,  the  transport inequality (definitions are recalled below) that one gets by standard techniques is as follows:  for every $f\ge 0$ on $M$ with $\int f\, d\sigma = 1$,
\begin{equation}\label{badtransport}
\mathcal W_{c}(f\, d\sigma, \sigma) \le \int f \log f \, d\sigma
\end{equation}
for the cost $c(d):=(n-1)d^2/2$ . Linearization of this inequality gives only a weak form of~\eqref{poinc} with $1/(n-1)$ in place  of the correct $1/n$. Let us note that by an abstract result of Otto and Villani~\cite{OV}, the  log-Sobolev inequality mentioned above  with the sharp constant $n$ implies that the transport inequality~\eqref{badtransport} holds with the cost $c(d):=n\, d^2/2$. As for the log-Sobolev inequality, it is not known how to reach this inequality using mass transport.

The difficulty is to properly quantify the interplay between dimension and non-zero curvature in the mass transportation techniques.

This was partly overcome in the work of Lott and Villani~\cite{LV07} (see also Chapters~20 and~21 of~\cite{V09}). There, the authors manage to prove some Sobolev-like inequalities under the so called "curvature-dimension condition $CD(K,n)$" that imply, after linearization,  sharp spectral bounds. To be precise, their assumption is that the metric measured space $(M,d,\sigma)$ satisfies a curvature-dimension lower bound which is defined in terms of uniform convexity  along optimal transport of a class of entropy functionals. From this assumption, they deduce a (not very natural)  Sobolev-like inequality. This inequality has no reason to be sharp when the curvature is nonzero, but after linearization it gives the correct Poincar\'e inequality~\eqref{poinc} (so in a sense it is sharp at first order).  Of course, it is known, by the properties of optimal transport  on manifolds (McCann's map), that a Riemannian manifold with condition~\eqref{riccibound} satisfies the curvature-dimension criterion.  So putting all together, we see that Lott and Villani's work is already an answer to the question on how to use  mass transport to derive some sharp dimensional inequalities. But of course, it is rather indirect, and no standard inequality that one could prove using optimal transport on a manifold is easy to extract from it. Actually, this is somehow the content of the "Open Problem 21.11" in Villani's book~\cite{V09}.

Our original motivation was to provide, in the particular case of a manifold, a different,  more  direct, approach  based on the geometric properties of  McCann's transport map. The aim was to find an inequality that contained the sharp bound~\eqref{poinc}. Eventually, we managed to establish a new, suitable, transport inequality, that is an inequality between an entropy functional and a transportation cost functional (we recommend the survey~\cite{GL} for background on transport inequalities).  The question of obtaining the sharp log-Sobolev inequality using mass transport remains.

Let us introduce the following classical dimensional entropy: given a probability density $f$ on $M$, meaning a Borel nonnegative function on $M$ with $\int f \, d\sigma = 1$, we put 
$$H_{n, \sigma}(f):= n \int \big( f- f^{1-1/n}\big)\, d\sigma = n - n\int f^{1-1/n} \, d\sigma.$$
Note that $H_{n, \sigma}$ is a nonnegative convex functional of $f$.

We will consider transportation costs given by  functions of the distance $d$ on $M$. Given a function $c:\R\to \R^+$ (or rather $c:[0,\pi]\to \R^+$ in our case),   the associated Kantorovich transportation cost between two Borel probability measures $\mu$ and $\nu$ on $M$ is defined by
$$\mathcal W_c(\mu, \nu):= \inf_\pi \iint  c\big(d(x,y)\big)\, d\pi(x,y)$$
where the infimum is taken over all probability measures $\pi$ on $M\times M$ projecting on $\mu$ and $\nu$, respectively.

In the proof of the Theorem below, will use McCann's map, which arises from an optimizer in the functional $\mathcal W_c$ when $c$ is the quadratic cost, $c(d)=d^2/2$; we shall recall McCann's result in detail later. However, let us emphasize that, although we will use this quadratic-optimal map, the cost in our transport inequality will be a different function of the distance. 

\def\S{{\rm S}}

Our cost function is defined for $d\in [0,\pi)$ by
$${c}_n(d):= n- \frac{\sin^{n-1}(d)}{\S_n(d)^{n-1}}- (n-1)\frac{\S_n(d)}{\tan(d)}$$
and at the limit by $c_n(\pi)=+\infty$, where $\S_n$ is the familiar function defined for $d\in [0, \pi]$ by
$$ {\rm S}_n(d):=  \Big( n\int_0^d\sin^{n-1}(s)  \, ds\, \Big)^{1/n}.$$
We have, as expected,  $c_n(0)=0$ (since $\S_n(t)\sim t$ at $0$) and $c_n(d)>0$ for $d>0$.

We can now state the transport inequality satisfied by the uniform measure $\sigma$ on  $M$.

\begin{theo*}
Let $M$ be an $n$-dimensional Riemannian manifold with positive Ricci curvature satisfying~\eqref{riccibound} and let  $\sigma$ be  its normalized Riemannian volume. 
Then, for every probability density $f$ on $M$ we have
$$\mathcal W_{c_n}(f\, d\sigma, \sigma) \le H_{n, \sigma} (f).$$
\end{theo*}
\medskip

We will see  that for small distances the cost $c_n(d(x,y))$ behaves  like $(n-1) d(x,y)^2/2$, so it may seem that we are back to the bad situation~\eqref{badtransport} where we were stuck with the constant $(n-1)$. However, the entropy $H_{n, \sigma}$ is~\emph{better}, i.e. smaller, than the usual entropy $\int f\log f\, d\sigma$ (note that $H_{n, \sigma}(f)\nearrow \int f\log(f) \, d\sigma$ as $n\to +\infty$) , and as a matter of fact we will reproduce the sharp Poincar\'e inequality. So there is an interesting trade-off between the cost and the entropy.
Incidentally, both sides of our inequality are zero when $n=1$ (which is a good sign), meaning that we don't derive any result on the torus $S^1$, although it might be possible,  by looking at first orders when $n\to 1$ and analyzing the proof below,  to guess what one should get in this case. 

The next section contains the proof of the Theorem. In the last section we  give some properties of the cost $c_n$ and we explain how to derive the sharp spectral gap inequality~\eqref{poinc} from the Theorem.

\section{Proof of the Theorem}

We start by recalling the result of McCann~\cite{M99}. Given two (compactly supported) probability densities $f$ and $g$ on a manifold $M$ with respect to $d\vol$, the Riemannian volume, there exists a Lipschitz function $\theta:M\to \R$ such that $-\theta$ is \emph{$c$-concave} and the map
$$T(x)=\exp_x(\nabla\theta(x))$$
pushes forward $f\, d\vol$ to $g\, d\vol$. The latter means that for every (bounded or nonnegative) Borel function $u$ on $M$,
$$\int u(y)\, g(y)\, d\vol(y) = \int u(T(x)) f(x)\, d\vol(x) .$$
The $c$-concavity of $-\theta$ is defined by the property that there exists a Lipschitz function $\psi$ such that $-\te(x) = \inf_y \{\psi(y) + d(x,y)^2/2\}$.  This implies (and is formally equivalent to) that at every point $x$ where $\theta$ is differentiable, and thus $y:=T(x)$ is uniquely defined, the function $v\to  \theta(v) + \frac12 d(v,y)^2 -\frac12 d(x,y)^2 $ achieves its minimum at $v=x$.  

Following a classical approach, the map $T$ is constructed by establishing that $\pi=(\id \times T) f\, d\vol$ is the optimizer for $\mathcal W_c(f\, d\vol, g\, d\vol)$ when $c$ is the quadratic cost. We will not use this  property, though.

As explained in~\cite{CMS01, CMS06}, it is possible to do, in a weak sense,  the change of variable $y=T(x)$ and to establish a pointwise Jacobian change of variable equation. To be precise, let us set, whenever it makes sense,  
$$dT_x := Y (H+\Hess_x \theta)$$
where, for fixed $x\in M$, the linear operators $Y:T_xM\to T_{T(x)}M$ and $H:T_xM\to T_xM$ are defined by
$$Y := Y_x=:d(\exp_x )_{\nabla\te(x)} \AND  H := H_x:= \Hess_x d^2 _{T(x)} /2,$$
 with the notation $d_y(\cdot)=d(y,\cdot)$ for fixed $y\in M$. Then, one has
$$f(x)= g(T(x))\, \det dT_x  \qquad f\, d\vol-a.e.$$ 
The set of points where this equation holds is contained in the set of $x\in M$  where $\te$ is differentiable at $x$ with  $\gamma(t):=\exp_x(t\nabla \te(x))$ being the unique minimizing geodesic between $x=\gamma(0)$ and $T(x)=\gamma(1)\notin {\rm cut}(x)$, and such that $\Hess_x \theta$ exists, in the sense of Aleksandrov for the Lipschitz (and locally semi-convex) function $\te$; later we shall use  that $\tr \Hess \theta=:\Delta \te \le \Delta_{\mathcal D} \theta$  where $\Delta_{\mathcal D} \theta $ is the distributional Laplacian of the Lipschitz function $\theta$. 
The $c$-concavity of $-\theta$ then implies the following, crucial monotonicity property of $T$, which holds  $f\, d\vol-a.e.$: 
\begin{equation}\label{monotony}
H+\Hess \theta \ge 0 
\end{equation}
In Euclidean space, $H=\id$ and we recover that $T(x)=x+\nabla \theta$ is the gradient of the convex function $|x|^2/2+\theta(x)$ (the Brenier map).

We refer the interested (or worried) reader to~\cite{CMS01, CMS06} where these facts are carefully stated and proved.

So, under the assumptions  of the Theorem, let $T(x) = \exp_x (\nabla \theta)$ be the McCann map  pushing $\sigma$ forward to $f\, d\sigma$. Denote the displacement distance by
$$\alpha(x) := d(x, T(x))= |\nabla \theta(x)|\in [0, \pi ].$$ 
The Jacobian equation satisfied almost everywhere is then
\begin{equation}\label{jacobian}f(T(x))^{-1}=\det\big[Y (H+\Hess_x \theta)\big]
\end{equation}
with $Y := Y_x=:d(\exp_x )_{\nabla\te(x)}$ and $H := H_x:= \Hess_x d^2 _{T(x)} /2$.

For $x\in M$ a point where equation~\eqref{jacobian} holds,  let $E_1 := \nabla\theta /|\nabla\theta|$ be the direction of transport, completed by $E_2, \ldots , E_n$ in order to have an orthonormal frame.
In this basis, the symmetric operator $H$ takes the form
$$ \begin{pmatrix}
  1  &   0  \\
 \smallskip
 0 &   K
 \end{pmatrix} $$
 and the classical Bishop comparison estimates (see e.g.~\cite{Petersen}) ensure that under~\eqref{riccibound} we have
$$\det Y \le \left( \frac{\sin ( \alpha)}{\alpha} \right)^{n-1}=:v_n(\alpha)^n\, \AND \quad
\mbox{tr}\, K \le (n-1) \frac{ \alpha }{ \tan( \alpha)}=: w_n(\alpha)     .$$
Of course, these inequalities are equalities when $M=S^n$, a case where $Y$ and $K$ can be computed explicitly (see~\cite{CE99}).

If we write $\Hess_x \theta=  \begin{pmatrix}
  a  &   b^t  \\
 \smallskip
 b &   M
 \end{pmatrix} $, where $M$ is a symmetric $(n-1)\times (n-1)$ matrix and $a:= \Hess_x \theta (E_1)\cdot E_1$ (all the quantities depend on $x$, of course), then we have
\begin{eqnarray*}
f(T(x))^{-1} &=&\det\Big[ Y\begin{pmatrix}
  1+a  &   b^t  \\
 \smallskip
 b &   K+M
 \end{pmatrix} \Big] 
 \le v_n(\alpha)^n  \det\begin{pmatrix}
  1+a  &   b^t  \\
 \smallskip
 b &   K+M
 \end{pmatrix}  
 \\
 & \le &  v_n(\alpha)^n \det\begin{pmatrix}
  1+a  &   0  \\
 \smallskip
 0 &   K+M
 \end{pmatrix} 
 \\
 & = &  v_n(\alpha)^n \det\begin{pmatrix}
  (1+a)\mu(\alpha)^{-(n-1)} &   0  \\
 \smallskip
 0 &   \mu(\alpha)K+\mu(\alpha)M
 \end{pmatrix} 
 \end{eqnarray*}
where $\mu$ is a numerical $C^1$ positive function defined on $[0, \pi]$ that will be fixed later. Note that $1+a\ge 0$ and $K+M\ge 0$ by~\eqref{monotony}.
Using  the arithmetic-geometric inequality, namely ${\det }^{1/n} \le$tr$/n$ on nonnegative matrices,  we then get that 
$$n \, f(T(x))^{-1/n}  
\le v_n(\alpha) \Big[ (1+a)\mu(\alpha)^{-(n-1)}+ \mu(\alpha)w_n(\alpha) + \mu(\alpha)(\Delta\te-a)\Big] $$
We integrate this inequality with respect to $\sigma$. Integration by parts gives
$$\int v_n(\alpha) \mu(\alpha) \Delta \theta \, d\sigma \le - \int (v_n \mu)'(\alpha) \nabla \alpha \cdot \nabla \te \, d\sigma.$$
When $\te$ is smooth, the previous equation is an equality; but as we explained above, the Laplacian we used is smaller than the distributional Laplacian, in general.

By construction,  $\nabla\alpha \cdot \nabla\te =\alpha  \Hess \te (E_1) \cdot E_1 = \alpha a$ (the fact that this property should be used to improve mass transportation techniques on manifold was suggested to us some years ago by Michael Schmuckenschl\"ager~\cite{Mi}).  So we find
\begin{eqnarray*}
n \int f^{1-1/n}\, d\sigma
&\le  &\int  \Big[v_n(\alpha) \mu(\alpha)^{-(n-1)} - \mu(\alpha) v_n(\alpha) -\alpha\cdot (v_n \mu)'(\alpha) \Big]\, a \,  d\sigma\\
& & \quad  + \int \big( \mu(\alpha)^{-(n-1)} + \mu(\alpha) w_n(\alpha)\big) v_n(\alpha)\, d\sigma.
\end{eqnarray*}
We now want to choose the numerical function  $ \mu$ such that for all $t\in [0, \pi)$,
\begin{equation}\label{ODEcond}
v_n(t) \mu(t)^{-(n-1)} - \mu(t) v_n(t) -t\, (v_n \mu)'(t) =0.
\end{equation}
Setting $h(t):= t \mu(t) v_n(t)$, the previous equation rewrites as
$$h'(t) = v_n(t) (h(t)/tv_n(t))^{-(n-1)}= v_n(t)^{n} t^{n-1} h(t)^{-(n-1)},$$
or equivalently
$$\frac1n (h^n)'(t) = \sin^{n-1}(t),$$
which suggests the choice $h=\S_n$.
So the function \emph{defined} by $\mu(t):= {\rm S}_n(t) /tv_n(t)$ satisfies~\eqref{ODEcond}, and consequently we have the desired inequality:
$$n \int f^{1-1/n}\, d\sigma
\le   \int \Big[ \frac{\sin^{n-1}(\alpha(x))}{\S_n(\alpha(x))^{n-1}}+ (n-1)\frac{\S_n(\alpha(x))}{\tan(\alpha(x))}\Big]\, d\sigma(x) .$$

\section{Further remarks}

We start with some properties of the function $c_n(\alpha)= n- \frac{\sin^{n-1}(\alpha)}{\S_n(\alpha)^{n-1}}- (n-1)\frac{\S_n(\alpha)}{\tan(\alpha)}$, $\alpha \in [0, \pi)$. 

First, observe that for $\alpha\in [0,\pi]$,
$\int_0^\alpha  \sin(s)^{n-1} \cos(s)\, ds \le \int_0^\alpha  \sin(s)^{n-1} \, ds \le \int_0^\alpha  s^{n-1}  \, ds$ so that 
$$\sin(\alpha) \le \S_n(\alpha) \le \alpha.$$
This implies that $c_n\ge 0$. It also gives that $0\le \frac{\alpha - \S_n(\alpha)}{\alpha^2} \le \frac{\alpha - \sin(\alpha)}{\alpha^2}$ and consequently, for $\alpha\to 0$, 
$$\S_n(\alpha) = \alpha +o(\alpha^2).$$
In turn, this gives the behavior of $c_n(\alpha)$ when $\alpha \to 0$:
\be\label{equic}
c_n(\alpha) \sim (n-1) \alpha^2/2.
\ee
To perform this series expansion of $c_n$,  write $\S_n(\alpha) = \alpha + a \alpha^3 + o(\alpha^3)$; the coefficient $a$ indeed disappears in the second order.
We believe (from numerical examples) that the function $c_n$ is convex on $[0,\pi]$. But since we don't need this property (which seems a bit more technical), we leave this question  for another time.

It is well known that the property~\eqref{equic}  of the cost is sufficient to derive by linearization,  from the corresponding transport inequality, a Poincar\'e type inequality. The standard procedure is to first state an infimal convolution inequality (for the Hamilton-Jacobi semi-group), obtained by dualizing the transportation cost and the entropy, and then to linearize  (see~\cite{GL}). Actually, it is enough to dualize only the transportation cost (we don't want to dualize the entropy, since eventually we will linearize it).

Recall the classical Kantorovich duality:  for two probability measures $\mu$ and $\nu$ on $M$ and for a cost $c$,
$$ \mathcal W_c(\mu,\nu)= \sup_\varphi \Big\{ \int Q_c(\varphi) d\mu - \int \varphi \, d\nu \Big\}$$
where the supremum is taken over all (Lipschitz) functions $\varphi:M\to \R$, and
$$\forall x \in M, \quad Q_c(\varphi)(x) := \inf_{y\in M} \big\{ \varphi(y) + c(d(x,y))\big\}.$$
Note that $Q_c(\varphi)\le \varphi$ (provided $c\ge 0$ and $c(0)=0$) and that the bigger the cost is in terms of $d(x,y)$, the closer $Q_c(\varphi)$ is to $\varphi$.

Let $g$ be a smooth function on $M$ with $\int g \, d\sigma = 0$, and $\eps>0$ small. Applying our transport inequality to the probability density $f=1+\eps \lambda g$ where $\lambda>0$ is a constant to be fixed later, and using the above-mentioned duality with the test function $\varphi= \eps g$ we get
\be\label{dualdel}
\int Q_{c_n} (\eps g) (1+\eps \lambda g) d\sigma - \int (\eps g) \, d\sigma \le H_{n,\sigma}(1+\eps \lambda g).
\ee
On one hand we have, for the entropy term,  uniformly on $M$,
$$  n\big[ (1+\eps \lambda g) - (1+\eps \lambda g)^{1-1/n}\big]= \eps \lambda g + \eps^2 \frac{n-1}{2n}(\lambda g)^2 + o(\eps^2).$$
On the other hand, because of~\eqref{equic} we have
$$Q_{c_n} (\eps g)  = \eps\big[ g - \eps\frac{1}{2(n-1)} |\nabla g |^2+ o(\eps)\big].$$
Putting these two expansions in~\eqref{dualdel}, we see that the orders $0$ and $1$ vanish (they have to, since the constant function $\mathbf 1$ is an equality case in the transport inequality), and the inequality between  the second orders reads as
$$\Big(\lambda - \frac{n-1}{2n} \lambda^2 \Big) \int g^2 \, d\sigma\le\frac{1}{2(n-1)} \int |\nabla g |^2\, d\sigma.$$
Picking $\lambda = \frac{n}{n-1}$ we get the sharp Poincar\'e inequality
$ \int g^2 \, d\sigma\le \frac{1}{n}\int |\nabla g |^2\, d\sigma$.

\medskip

\noindent Dario Cordero\,-Erausquin: \verb?cordero@math.jussieu.fr?\\
Institut de Mathématiques de Jussieu and Institut Universitaire de France, \\
Universit\'e Pierre et Marie Curie -- Paris 6, \\
4, place Jussieu, 
F-75252 Paris Cedex 05,
France

\end{document}